\documentclass[11pt,a4paper]{article}
\begin{document}

\title{History of the formulas and algorithms for $\pi$}

\author{Jes\'{u}s Guillera Goyanes}
\date{}

\maketitle

\newtheorem{conj}{Conjecture}[section]

\begin{abstract}
Throughout more than two millennia many formulas have been obtained, some of them beautiful, to calculate the number $\pi$. Among them, we can find series, infinite products, expansions as continued fractions and expansions using radicals. Some expressions which are (amazingly) related to $\pi$ have been evaluated. In addition, a continual battle has been waged just to break the records computing digits of this number; records have been set using rapidly converging series, ultra fast algorithms and really surprising ones, calculating isolated digits. The development of powerful computers has played a fundamental role in these achievements of calculus.
\end{abstract}

\section{First formula: Archimedes' algorithm}

For a period of approximately 1800 years, Archimedes' algorithm, proved around 250 BC, was the most efficient way to calculate $\pi$. The idea consisted of considering a circle with unit diameter with regular circumscribed and inscribed polygons of $3 \cdot 2^n$ sides. Let $a_n$ and $b_n$ denote their perimeters, respectively. He proved the following relations using only geometrical reasoning:
\begin{equation}\label{archimedes}
a_1=2 \sqrt{3}, \quad b_1=3, \qquad a_{n+1}=\frac{2 a_n
b_n}{a_n+b_n}, \qquad b_{n+1}=\sqrt{a_{n+1} b_n}.
\end{equation}
Obviously, $b_n < \pi < a_n$ and both sequences $a_n$ and $b_n$ converge to $\pi$. It is an algorithm which nowadays can be easily proved using elementary trigonometry. Indeed, if $k_n=3 \cdot 2^n$, we can write
\[ a_n=k_n \tan \frac{\pi}{k_n}, \qquad b_n=k_n \sin \frac{\pi}{k_n}. \]
Then
\[ \frac{2 a_ n b_ n}{a_n+b_n}=2k_n \frac{\tan \frac{\pi}{k_n} \sin \frac{\pi}{k_n}}{\tan \frac{\pi}{k_n}+
\sin \frac{\pi}{k_n}}= 2k_n \tan \frac{\pi}{2k_n}=a_{n+1}. \]
On the other hand,
\[ \sqrt{a_{n+1} b_n}=\sqrt{2k_n \tan \frac{\pi}{2k_n} \cdot k_n \sin \frac{\pi}{k_n}}=
2k_n \; \sqrt{\tan \frac{\pi}{2k_n} \sin \frac{\pi}{2k_n} \cos
\frac{\pi}{2k_n}}=b_{n+1}. \]

\par Convergence in this algorithm is linear; after $5$ steps we get $3$ digits of  $\pi$. To see this, we use the inequalities  $\tan x < x$ and $1-\cos x <
\frac{x^2}{2}$
to get
\[
a_n-b_n=k_n \tan \frac{\pi}{k_n} \left( 1-\cos \frac{\pi}{k_n} \right)< \frac{\pi^3}{2 k_n^2}=
\frac{\pi^3}{18 \cdot 2^{2n}}<\frac{\pi^3}{18\cdot 1000^{n/5}}.
\]
With $n=7$, this algorithm gives the estimate $3.1415< \pi< 3.1417$.

\par It is enlightening to point out that during the age of Archimedes neither decimal notation nor any other positional notation was known, so he stated the result he obtained in terms of polygons with $96$ sides ($n=5$) using fractions:
\[ 3 + \frac{10}{71} < \pi < 3 + \frac{1}{7}. \]
This approximation determines two-place accuracy; this means $\pi\approx3.14$.

\section{No great news until 1800 years later}

For a period of approximately $1800$ years all formulas that were obtained for $\pi$ were based on the idea of approximating a circle using regular polygons, and they took the shape of iterative algorithms. In $1593$ Vi\`{e}te obtained a different expression by means of a pretty expansion with infinitely many nested square roots:
\begin{equation}\label{viete}
\frac{2}{\pi}= \sqrt{\frac{1}{2}} \cdot
\sqrt{\frac{1}{2}+\frac{1}{2} \sqrt{\frac{1}{2}}} \cdot
\sqrt{\frac{1}{2}+\frac{1}{2} \sqrt{\frac{1}{2}+ \frac{1}{2}
\sqrt{\frac{1}{2}}}} \cdots.
\end{equation}
As with all the previous methods, Vi\`{e}te's is based on geometry and consists of considering a circle whose diameter is $2$ and approximating its area, that is, $\pi$, by the area of an inscribed square, the area of an inscribed octagon, the area of an inscribed 16-sided polygon, and so on. In this way Vi\`{e}te found $\pi$ accurately to nine places: $3.141592653$.

\par Around 1650, analytical methods of Descartes and Fermat, who were forerunners of the development of Calculus which would be carried out soon, started replacing geometrical methods. One of the formulas for $\pi$ at this stage is the famous infinite product due to Wallis (1655):
\begin{equation}\label{wallis}
\frac{2}{\pi}=\frac{1}{2} \cdot \frac{3}{2} \cdot \frac{3}{4} \cdot \frac{5}{4} \cdot \frac{5}{6} \cdot \frac{7}{6} \cdot \frac{7}{8} \cdot \frac{9}{8} \cdots.
\end{equation}

\par In 1999, Osler proved a formula which includes as particular cases both Vi\`{e}te's formula and Wallis' infinite product. We will see this in section $10$.

\section{The early stages of the Calculus: New Formulas for $\pi$}

Between 1665 and 1680, Newton and Leibniz independently developed the first basic notions of Calculus, which supplied us with a powerful tool not only for mathematics but also for physics and other sciences. These methods, as it could not be otherwise, also had repercussions in obtaining new formulas for $\pi$. For example, one obtained by Newton is
\begin{equation}\label{newton}
\pi= \frac{3 \sqrt{3}}{4}+2-\frac{3}{4} \sum_{n=0}^{\infty}
\frac{{2n \choose n}}{(n+1)(2n+5) 16^n},
\end{equation}
by which he approximated $\pi$ to $15$ places in 1665. Newton reached this formula by noticing that the area which is bounded by a circle with the equation
\[ \left(x-\frac{1}{2} \right)^2+y^2=\frac{1}{4} \]
and the straight lines $x=\frac{1}{4}$ and $y=0$ is equal to the difference between the area of a circular sector that takes up an angle of $60^{o}$ and its corresponding triangle. Therefore
\[ \int_0^{1/4} y \; dx = \int_0^{1/4} \sqrt{x-x^2} dx=
\int_0^{1/4} x^{1/2} (1-x)^{1/2} dx=\frac{\pi}{24}-\frac{\sqrt
3}{32}. \]
Newton's formula is obtained from here by expanding $(1-x)^{1/2}$ using the binomial formula, which he had proved, and then integrating term by term. Formulas such as (\ref{newton}) are nowadays known as binomial sums.

\par A series similar to Newton's one but simpler is
\begin{equation}
\frac{\pi}{3}=\sum_{n=1}^{\infty} \frac{{2n \choose n}}{(2n+1)
16^n},
\end{equation}
which we get by substituting $x=\frac{1}{4}$ into the series expansion (see \cite{borwein1}, p. 386)
\[ \frac{\arcsin 2x}{2x}=\sum_{n=1}^{\infty} \frac{{2n \choose n}}{2n+1} x^{2n}. \]

\par Another interesting example of a binomial sum is the following obtained by Comtet in 1974:
\begin{equation}
\pi^4= \frac{3240}{17} \sum_{n=1}^{\infty} \frac{1}{n^4 {2n \choose
n}}
\end{equation}
which can be deduced from the identity \cite{borwein3}
\[
\sum_{n=1}^{\infty} \frac{1}{n^4 {2n \choose n}}=2 \int_0^{\pi/3} \theta \log^2 \left( 2 \sin \frac{\theta}{2}
\right) d \theta.
\]

\par In 1671 Gregory integrating the following geometric series term by term:
\[ 1-x^2+x^4-x^6+\cdots=\frac{1}{1+x^2}, \qquad |x|<1, \]
got the formula (which has his name)
\[
\arctan x= x-\frac{x^3}{3}+\frac{x^5}{5}-\frac{x^7}{7}+\cdots,
\qquad |x| \leq 1,
\]
which would become known as the power series expansion of the function $\arctan x$. Substituting $x=1$, Leibniz got in 1674 the series for $\pi$

\begin{equation}\label{leibniz}
\frac{\pi}{4}= 1-\frac{1}{3}+\frac{1}{5}-\frac{1}{7}+\cdots.
\end{equation}
It is strange that Gregory did not do it himself, probably because he was not interested in slowly convergent series. Sharp put $x=\sqrt{3}/3$ in Gregory's formula and obtained the series
\begin{equation}\label{sharp}
\frac{\pi}{6}=\frac{\sqrt{3}}{3} \left( 1-\frac{1}{3 \cdot 3}+
\frac{1}{5 \cdot 3^2}-\frac{1}{7 \cdot 3^3}+\cdots \right)
\end{equation}
which allowed him to achieve the record of 71 digits for $\pi$ in 1699. It is surprising that Gregory's formula and Sharp's particular case have been found in the Sanskrit texts that the Indian learned person Nilakantha Somayaji had written 200 years before.

\par In 1665 Brouncker got the pretty expansion as a continued fraction
\begin{equation}
\frac{4}{\pi}= 1 + {\displaystyle \frac{1^2}{2+ {\displaystyle
\frac{3^2}{2+ {\displaystyle \frac{5^2}{2+ {\displaystyle
\frac{7^2}{2+ {\displaystyle \frac{9^2}{2+\ddots}}}}}}}}}}.
\end{equation}
About $100$ years later, Euler proved it in a clever way. He starts by proving the equality
\[
a_1 + a_1 a_2 + a_1 a_2a_3+ a_1 a_2 a_3 a_4 \cdots
=\frac{a_1}{1-{\displaystyle \frac{a_2}{1+a_2-{\displaystyle
\frac{a_3}{1+a_3-{\displaystyle \frac{a_4}{1+a_4- \ddots}}}}}}},
\]
which motivates him to rewrite Gregory's series as
\[ \arctan x=x+x \left(\frac{-x^2}{3}\right)+
x \left(\frac{-x^2}{3}\right) \left( \frac{-3x^2}{5} \right)+ x
\left(\frac{-x^2}{3}\right) \left( \frac{-3x^2}{5} \right) \left(
\frac{-5x^2}{7} \right) + \cdots
\]
and then he gets Brouncker's formula by substituting $x=1$ and using the continued fraction.

\par A recent example of expanding $\pi$ as a continued fraction is
\begin{equation}\label{lange}
\pi= 3 + {\displaystyle \frac{1^2}{6+ {\displaystyle \frac{3^2}{6+
{\displaystyle \frac{5^2}{6+ {\displaystyle \frac{7^2}{6+
{\displaystyle \frac{9^2}{6+\ddots}}}}}}}}}},
\end{equation}
obtained by Lange \cite{lange} in 1999.

\section{Machin type series}

In 1706, John Machin proved the following formula (which is known today as Machin's formula):
\begin{equation}\label{machin1}
\frac{\pi}{4}=4 \arctan \frac{1}{5}-\arctan \frac{1}{239}.
\end{equation}
Using Gregory's formula he rewrote it as a series and in this way he got the first 100 digits of $\pi$. The following identities are similar, but they lead to more slowly convergent series:
\begin{equation}\label{tipomachin}
\frac{\pi}{4}=\arctan \frac{1}{2}+ \arctan \frac{1}{3}, \qquad
\frac{\pi}{4}=2 \arctan \frac{1}{3}+ \arctan \frac{1}{7}.
\end{equation}

\par The proof of a known Machin type series expression consists of an ordinary checking of an identity between elementary operations with complex numbers. We use
\[ \frac{(5+i)^4}{239+i}=2+2i \]
to prove Machin's formula; it is enough to make the arguments equal.

\par Gauss (1777-1855) also took part in the search for rapidly convergent series for $\pi$. He proved the following equality:
\begin{equation}\label{gaussmachin}
\frac{\pi}{4}= 12 \arctan \frac{1}{18} +8 \arctan \frac{1}{57} -5
\arctan \frac{1}{239}.
\end{equation}
It is like the preceding ones, but it involves one more term.
\par In $2002$ Kanada and his team, using a certain type of trick to optimize computer calculations, got a record computing $\pi$ by means of the series
\begin{equation}\label{takano}
\frac{\pi}{4}= 12 \arctan \frac{1}{49}+32 \arctan \frac{1}{57}-5
\arctan \frac{1}{239}+12 \arctan \frac{1}{110443},
\end{equation}
\begin{equation}\label{stormer}
\frac{\pi}{4}= 44 \arctan \frac{1}{57}+7 \arctan \frac{1}{239}-12
\arctan \frac{1}{682}+24 \arctan \frac{1}{12943}.
\end{equation}
The first one was obtained by Takano in 1982 (he is a high school teacher and a song composer). The second one was obtained by St\"{o}rmer in $1896$. Kanada's record consisted of more than a billion of digits of $\pi.$

\section{Formulas by the Swiss mathematician L. Euler}

Great mathematicians, Leibniz and the Bernoulli brothers Jacob (1654-1705) and Johann (1667-1748) among them, had attempted to evaluate (with no success) the infinite sum
\begin{equation}\label{euler}
\frac{1}{1^2}+\frac{1}{2^2}+\frac{1}{3^2}+\frac{1}{4^2}+\frac{1}{5^2}+\cdots.
\end{equation}
Finally, Euler proved in 1736 that its value is $\pi^2/6.$ To do it he started with the series
\[ \frac{\sin x}{x}=1-\frac{x^2}{3!}+\frac{x^4}{5!}-\frac{x^6}{7!}+\cdots \]
and the equation
\[ 1-\frac{x^2}{3!}+\frac{x^4}{5!}-\frac{x^6}{7!}+\cdots = 0. \]
The solutions of this equation should be the solutions of $\sin x=0$, except $x=0$, that is: $\pm \pi$, $\pm 2\pi$, $\pm 3 \pi$, $\ldots$. He imagined the infinite series as a polynomial and he broke it down into factors
\[
1-\frac{x^2}{3!}+\frac{x^4}{5!}-\frac{x^6}{7!}+\cdots =
\left( 1-\frac {x^2}{\pi^2} \right) \left( 1-\frac {x^2}{4 \pi^2}
\right) \left( 1-\frac {x^2}{9\pi^2} \right) \cdots.
\]
He finally got the proof by setting the $x^2$ coefficients equal. He also proved many similar formulas. Among them we can find
\begin{equation}\label{eulerotras}
\sum_{n=1}^{\infty} \frac{1}{(2n+1)^2}=\frac{\pi^2}{8}, \qquad
\sum_{n=1}^{\infty} \frac{1}{n^4}=\frac{\pi^4}{90}, \qquad
\sum_{n=0}^{\infty} \frac{(-1)^{n}}{(2n+1)^3}=\frac{\pi^3}{32}
\end{equation}
and the not-so-similar
\begin{equation}\label{eulerotramas}
\sum_{n=1}^{\infty} \frac{ {\displaystyle
1+\frac{1}{2}+\frac{1}{3}+\cdots\frac{1}{n}}}{n^3}=\frac{\pi^4}{72}.
\end{equation}
Nowadays formulas like the last one are known as Euler sums \cite{flajolet}. The following one has an interesting history
\begin{equation}\label{sumaeuler}
\sum_{n=1}^{\infty} \frac{ {\displaystyle
\left(1+\frac{1}{2}+\frac{1}{3}+ \cdots \frac{1}{n}
\right)^2}}{n^2}= \frac{17 \pi^4}{360},
\end{equation}
because it was observed numerically by the student Au Yeung. J. Borwein, his professor, was skeptical of the result but after checking it with greater precision he got to prove it. Now, we know that the formula had been proved long before by H. Sandham \cite{otto}. However, the rediscovery of this forgotten formula has had an important effect because it has motivated the organized study of these kind of formulas.

\par While studying the divergence of the series of reciprocals of the prime numbers, Euler considered geometric progressions:
\[1+\frac{1}{p^2}+\frac{1}{p^4}+\frac{1}{p^6}+\frac{1}{p^8}+\cdots= \frac{1}{1-\frac{1}{p^2}}, \] and he observed that by multiplying over all primes $p$ we get
\begin{equation}\label{euler2}
\sum_{n=1}^{\infty} \frac{1}{n^2} =\frac{\pi^2}{6}=\prod_{p}
\left( 1-\frac{1}{p^2} \right)^{-1},
\end{equation}
which links the number $\pi$ with prime numbers.

\par From 1748 he initiated the study of the lemniscate and discovered the formula
\begin{equation}\label{prod-integrals}
\int_0^1 \frac{dt}{\sqrt{1-t^4}} \cdot \int_0^1 \frac{t^2dt}{\sqrt{1-t^4}}=\frac{\pi}{4},
\end{equation}
which establish a remarkable relationship between two elliptic integrals.

\section{An outstanding forgotten algorithm: C. F. Gauss' Formula}

When he was only 14 years old, Gauss considered the following iterative process:
\[ a_0=a, \quad b_0=b, \quad a_{n+1}=\frac{a_n+b_n}{2}, \qquad b_{n+1}= \sqrt{a_n b_n} \]
and he defined the arithmetic-geometric mean of $a$ and $b$ as the common limit
\[ M(a,b)=\lim_{n \to \infty} a_n=\lim_{n \to \infty} b_n. \]
In 1800, his investigations concerning the arithmetic-geometric mean lead him to the wonderful formula (see \cite{arndt}, p. 94-102)
\begin{equation}\label{gauss}
\sum_{k=0}^{\infty} 2^{k} (a_k^2-b_k^2)=1-\frac{2M^2}{\pi},
\end{equation} where $a_0=1$, $b_0=1/\sqrt 2$ and $M$ is the arithmetic-geometric mean of $a_0$ and $b_0$.

We sketch a proof of it using the elliptic integrals
\[
I_n= \int_0^{\pi/2} \frac{d \theta}{\sqrt{a_n^2 \cos^2 \theta+b_n^2\sin^2 \theta}}, \qquad
L_n= \int_0^{\pi/2} \frac{\sin^2 \theta d \theta}{\sqrt{a_n^2 \cos^2 \theta+b_n^2\sin^2 \theta}}.
\]
The first step consist in proving that $I_n=I_{n+1}$ by means of clever substitutions of the integration variable. This result has an important consequence, namely:
\[ I_0=I_n=\lim_{n \to \infty} I_n=\frac{\pi}{2M}. \]
The second step consist in deriving, again by clever substitutions, the relation
\[ c_n^2 L_n-2 c_{n+1}^2 L_{n+1}=\frac{1}{2} c_n^2 I_n, \quad {\rm where} \quad c_n^2=a_n^2-b_n^2. \]
Multiplying by $2^n$ and summing for $n \geq 0$, we obtain
\[ c_0^2 L_0=\frac{1}{2} I_0 \sum_{n=0}^{\infty} 2^n c_n^2. \]
The last step is the evaluation $I_0(I_0-L_0)=\pi/2$, which follows easily from Euler's formula (\ref{prod-integrals}).

\par For no clear reasons, Gauss did not make use of his formula to calculate $\pi$. In 1970, it was rediscovered independently by Eugene Salamin and Richard Brent \cite{salamin}, \cite{brent} who expressed it appropriately for iterative calculation:
\[ \lim_{n \to \infty} \frac{2a_n^2}{1- {\displaystyle \sum_{k=0}^n 2^k (a_k^2-b_k^2)}}=\pi. \]
It can be proved that the convergence to $\pi$ is quadratic, that is, from $n$ to $n+1$ the number of correct digits is doubled. Brent and Salamin programmed this algorithm on a computer. It let them reach in 1976 a record approximating $\pi$ to $3$ million digits, while the record until that moment had been $1$ million digits, achieved using Machin-type formulas.

\section{Series for $1/\pi$ from an incredible Indian mathematical genius: S.~Ramanujan}

The perspicacity and ability of Ramanujan allowed him to discover, among many other things, 17 extraordinary series for $1/\pi$ completely different from those known until that moment. These new series were published in 1914 in his celebrated paper \cite{ramanujan}. In it he considers the complete elliptic integral of the first kind
\[
\int_0^{\pi/2}\frac{d \phi}{\sqrt{1-k^2 \sin^2 \phi}}=1+\left( \frac{1}{2} \right)^2 k^2 + \left( \frac{1 \cdot 3}{2 \cdot 4} \right)^2 k^4 + \cdots.
\]
The number $k$ is called the modulus and $k'=\sqrt{1-k^2}$ is called the complementary modulus. He then gives his definition of what ordinary modular equations are: They express the algebraic relations which hold between $k$ and $l$ when, for a given rational $r$, we have $r L'/L=K'/K$, or $Q^r=q$, where $q=e^{-\pi K'/K}$, $Q=e^{-\pi L'/L}$ and $K$, $K'$, $L$ and $L'$ denote complete elliptic integrals with the moduli $k$, $k'$, $l$ and $l'$. Supposing $k=l'$ and $k'=l$, we have $q=e^{-\pi \sqrt{r}}$. Developing this theory \cite{ramanujan}, \cite{berndt} he proved the formulas
\begin{equation}\label{rama42n+5}
\sum_{n=0}^{\infty}  \frac{(2n)!^3}{n!^6} \frac{1}{2^{8n}} (6n+1)=\frac{4}{\pi}, \qquad
\sum_{n=0}^{\infty} \frac{(2n)!^3}{n!^6} \frac{1}{2^{12n}} (42n+5)=\frac{16}{\pi}
\end{equation}
and
\begin{equation}\label{rama-zalg}
\sum_{n=0}^{\infty} \frac{(2n)!^3}{n!^6} \left( \frac{3-\sqrt{5}}{16} \right)^{4n} ((42\sqrt{5}+30)n+(5\sqrt{5}-1))=\frac{32}{\pi}.
\end{equation}
He says that the other 14 examples he gives belong to three alternative theories: the quartic, the cubic and the sextic. He defines the corresponding functions $K$, but it is a pity that he does not develop these theories. One of his examples in the quartic theory is
\begin{equation}\label{ramanujan4}
\sum_{n=0}^{\infty} \frac{(4n)!}{n!^4} \frac{1}{396^{4n}}(26390n+1103)=\frac{9801 \sqrt{2}}{4 \pi}.
\end{equation}
In the cubic theory he gives, for example, the series
\begin{equation}\label{ramanujan5}
\sum_{n=0}^{\infty} \frac{(2n)!(3n)!}{n!^5} \frac{1}{1458^n}(15n+2)=\frac{27}{4 \pi}.
\end{equation}
And in the sextic
\begin{equation}\label{ramanujan6}
\sum_{n=0}^{\infty} \frac{(6n)!}{(3n)! n!^3} \frac{1}{54000^n}(11n+1)=\frac{5 \sqrt{15}}{6 \pi}.
\end{equation}
Series (\ref{ramanujan4}) is incredibly fast: each term contributes 8 digits of $\pi$. In $1985$, before a proof of this formula was known, W. Gosper programmed it on a computer and got a record of $17$ million digits of $\pi$.
\par The Borwein brothers, Jonathan and Peter, were the first to give complete and rigorous proofs, published in 1987, of the 17 Ramanujan series \cite{borwein1}. At about the same time the Chudnovsky brothers, David and Gregory, developed in \cite{chudnovsky} a theory to prove the following amazing series formula, which was not discovered by Ramanujan:
\begin{equation}\label{chudnovsky}
\sum_{n=0}^{\infty} \frac{(6n)!}{(3n)! n!^3}
\frac{(-1)^n}{640320^{3n}}
(545140134n+13591409)=\frac{\sqrt{640320^3}}{12 \pi}.
\end{equation}
Each term contributes $14$ digits of $\pi$ and they used it to beat records of the calculation of $\pi$ in $1989$, $1991$ and $1994$, the last one with 4044 million digits. In addition, they proved that there formula (\ref{chudnovsky}), which corresponds to $r=163$, has the fastest convergence rate among all those Ramanujan series that involve a rational raised to the power $n$. They arrived at this conclusion by noticing that $d=163$ is the largest positive integer such that the field $\mathbf{Q}(\sqrt{-d})$ has the unique factorization property \cite{chudnovsky}. In the year 2001, H. H. Chan, W. C. Liaw and V. Tan found identities \cite{chan1} which allowed them to prove the following formula, which also had not been discovered by Ramanujan:
\begin{equation}\label{chanhh}
\sum_{n=0}^{\infty} \frac{(3n)!(2n)!}{n!^5} \frac{(-1)^n}{300^{3n}}
(14151n+827)=\frac{1500 \sqrt{3}}{\pi}.
\end{equation}
Very recently  N. Baruah and B. Berndt in \cite{baruah1} and \cite{baruah2} following more closely the initial ideas of Ramanujan as presented in Section 13 of the celebrated paper \cite{ramanujan} got a much simpler treatment of the theory, and proved many new series of this type. An excellent survey on Ramanujan's series is \cite{baruah0}.

\section{The Borwein brothers' quartic algorithm}

The fourth root of the elliptic lambda function and the elliptic alpha function \cite{borwein1}, that is,
\[
s(q)=\sqrt[4]{\lambda(q)}=\frac{\theta_2(q)}{\theta_3(q)}, \quad
\alpha(q)=\frac{1}{\pi} \cdot \frac{1}{\theta_3^4(q)} \left[ 1+(4q \ln q)
\frac{\theta_4^{\,'}(q)}{\theta_4(q)} \right],
\]
where
\[
\theta_2(q)=\sum_{n=-\infty}^{n=\infty} q^{{(n+1/2)}^2}, \quad
\theta_3(q)=\sum_{n=-\infty}^{n=\infty} {q^{n^2}}, \quad
\theta_4(q)=\sum_{n=-\infty}^{n=\infty} (-1)^n{q^{n^2}}
\]
are the Jacobi elliptic theta functions, satisfy the following wonderful properties:
\[ s(q^4)=\frac{1-\sqrt[4]{1-s^4(q)}}{1+\sqrt[4]{1-s^4(q)}}, \]
\[ \alpha(q^4)=\alpha(q) [1+s(q^4)]^4+\frac{4 \ln q}{\pi} s(q^4) [1+s(q^4)+s^2(q^4)]. \]
These identities connect the function values at $q$ with the values at $q^4$, so they allow us to obtain an amazing Borwein family of quartic algorithms \cite{borwein2}. For example, if we let
\[ q=e^{- 2 \pi \cdot 4^n}, \quad s_n=s(q), \quad t_n=\alpha(q), \]
then the values of $s_0$ and $t_0$ are algebraic, and we get the following quartic algorithm:
\[ s_0=\sqrt{2}-1, \quad t_0=6-4 \sqrt{2}, \]
\begin{equation}\label{algocuartico}
s_{n+1}=\frac{1-\sqrt[4]{1-s_n^4}}{1+\sqrt[4]{1-s_n^4}},
\end{equation}
\[ t_{n+1}=t_n(1+s_{n+1})^4-2^{2n+3} s_{n+1} (1+s_{n+1}+s_{n+1}^2), \] where $t_n$ tends to $1/\pi$, and the number of accurate digits is quadrupled in each iteration. In the same way, if we let
\[ q=e^{-\pi \cdot 4^n}, \quad s_n=s(q), \quad t_n=\alpha(q), \]
we obtain a quartic algorithm with initial values $s_0=1/\sqrt[4]{2}$ and $t_0=1/2$.
A simple proof of this one is given in \cite{guillera6} using only Gauss' formula (\ref{gauss}) and elementary algebra. However, although this proof is much simpler than the Borweins', there method is more general. In addition, J. Borwein and P. Borwein have obtained other quadratic, cubic, etc algorithms. However, it seems that the quartic algorithms are the most efficient.
\par We might think that $\pi$ calculation records would always be broken using quadratic, quartic, etc. algorithms. However, that is not true, because such algorithms do not auto-correct their errors and they require doing all the operations to the final desired accuracy. Other ways have been found to improve the accuracy of calculation of series without affecting the final result, allowing them to beat records even though they have the drawback of being only linear.

\section{The jumbled calculation of $\pi$ digits}

In 1995 Peter Borwein and Simon Plouffe realized that the series
\[
\ln 2=\frac{1}{1 \cdot 2}+\frac{1}{2 \cdot 2^2}+\frac{1}{3 \cdot 2^3}+\frac{1}{4 \cdot 2^4}+\cdots=
\sum_{n=1}^{\infty} \frac{1}{n 2^n}
\]
could be used to compute one of the binary digits of $\ln 2$ without calculating the previous digits, by means of the following ingenious formula, where $\{ x \}$ means the fractional part of $x$:
\[
\left\{ 2^k \ln 2 \right\}=\left\{ \left\{ \sum_{n=0}^k \frac{2^{k-n}}{n} \right\}+ \sum_{n=k+1}^{\infty} \frac{1}{n 2^{n-k}} \right\}=\left\{ \left\{ \sum_{n=0}^k \frac{2^{k-n} \, {\rm mod} \, n}{n} \right\}+ \sum_{n=k+1}^{\infty} \frac{1}{n 2^{n-k}} \right\}.
\]
It lets us efficiently obtain the binary digit of $\ln 2$ in position $k+1$. The main remark to carry out all the involved calculation is that the numerator of the first sum, that is, $2^{k-n} \, {\rm mod} \, n,$ can be computed very quickly. For example, the calculation of $2^{65}$ does not require $65$ multiplications but only $7$ if we do it as $((((((2^2)^2)^2)^2)^2)^2) \cdot 2$ and a very efficient algorithm is used to obtain for example $2^{65} \, {\rm mod} \, 7$. It consists of reducing modulo $7$ after each multiplication instead of doing them in a row.
\par In 1989 the Borwein brothers claimed that obtaining just one binary digit of $\pi$ could not be easier than computing all the previous digits. However, after this discovery for $\ln 2$, it was not so clear and the team consisting of D. Bailey, P. Borwein and S. Plouffe began to search, using the PSLQ algorithm, for a series for $\pi$ that would allow them to calculate isolatedly its base $2$ digits (see \cite{delahaye}, chapter 8).
\par The PSLQ algorithm is a numerical algorithm which operates to a prearranged accuracy. It links an output of $n$ integers $a_1, a_2,\ldots,a_n,$ where $a_i\neq0$ for any $i$, with an input of $n$ real numbers $x_1,x_2,\ldots,x_n,$ and
\[ a_1 x_1 + a_2 x_2 + \cdots + a_n x_n=0 \]
is verified.

\par Success arrived taking
\[ x_j=\sum_{n=0}^{\infty} \frac{1}{(8n+j)16^n}, \qquad j=1,2, \ldots, 7, \] and trying to find, using PSLQ, a linear relation with integer coefficients using as input $\pi$, $x_1$, $x_2$, $x_3$, $x_4$, $x_5$, $x_6$, $x_7$. The output obtained was this sequence of integers: $-1$, $4$, $0$, $0$, $-2$, $-1$, $-1$, $0$. That is, they found a series \cite{bailey} which they called BBP type (BBP are the initials of the team members' names)
\begin{equation}\label{plouffe}
\pi = \sum_{n=0}^{\infty} \frac{1}{16^n} \left(
\frac{4}{8n+1}-\frac{2}{8n+4}-\frac{1}{8n+5}-\frac{1}{8n+6} \right)
\end{equation}
that lets us calculate $\pi$'s hexadecimal digits (hence also its binary ones) in isolation. Later, a proof of this formula was obtained from the easy result
\[ \int_0^{1/\sqrt{2}} \frac{x^{k-1}}{1-x^8} dx =
\int_0^{1/\sqrt{2}} \sum_{n=0}^{\infty} x^{k-1+8n} dx =
\frac{1}{2^{k/2}} \sum_{n=0}^{\infty} \frac{1}{16^n (8n+k)}.
\]
Namely, we have
\[ \sum_{n=0}^{\infty} \frac{1}{16^n} \left(\frac{4}{8n+1}-\frac{2}{8n+4}-\frac{1}{8n+5}-\frac{1}{8n+6} \right)=
\int_0^{1/\sqrt{2}} \frac{4 \sqrt{2}-8x^3-4 \sqrt{2}
x^4-8x^5}{1-x^8} dx,
\]
and this last integral can be easily evaluated by substituting $y=x \sqrt{2}$, as follows:
\[ 16 \int_0^1 \frac{4-2y^3-y^4-y^5}{16-y^8}dy=\int_0^1 \frac{16y-16}{(y^2-2)(y^2-2y+2)} dy=\pi. \]

\par As soon as they announced the amazing formula (\ref{plouffe}), other investigators wanted to find other beautiful BBP type formulas. The easiest one for $\pi$ was found and proved by Viktor Adamchik and Stan Wagon (see \cite{zhukov} pp. 66 and 67)
\begin{equation}\label{adamchik}
\pi= \sum_{n=0}^{\infty} \frac{(-1)^n}{4^n} \left(
\frac{2}{4n+1}+\frac{2}{4n+2}+\frac{1}{4n+3} \right).
\end{equation}
It can be proved similarly to the previous one, just taking into account the identity
\[
\int_0^{1/\sqrt{2}} \frac{x^{k-1}}{1+x^4} dx =
\int_0^{1/\sqrt{2}} \sum_{n=0}^{\infty} (-1)^n x^{k-1+4n} dx =
\frac{1}{2^{k/2}} \sum_{n=0}^{\infty} \frac{(-1)^n}{4^n (4n+k)}.
\]

\par In $1997$, Fabrice Bellard proved the formula
\begin{equation}\label{bellard}
\pi= 4 \sum_{n=0}^{\infty} \frac{(-1)^n}{4^n (2n+1)}- \frac{1}{64}
\sum_{n=0}^{\infty} \frac{(-1)^n}{1024^n} \left(
\frac{32}{4n+1}+\frac{8}{4n+2}+\frac{1}{4n+3} \right),
\end{equation}
and used it to brake a calculation record finding that the billionth binary digit of $\pi$ is $1$. His formula is now used in almost every computation of isolated binary digits of $\pi$ and the actual record belongs to C. Percival.

\par There are also known BBP formulas for other values related to $\pi$, for example,
\begin{equation}\label{bbpotra1}
\pi \sqrt{3}= \frac{9}{32} \sum_{n=0}^{\infty} \frac{1}{64^n} \left(
\frac{16}{6n+1}+\frac{8}{6n+2}-\frac{2}{6n+4}-\frac{1}{6n+5} \right)
\end{equation}
and
\begin{equation}\label{bbpotra2}
\pi^2= \frac{9}{8} \sum_{n=0}^{\infty} \frac{1}{64^n} \left(
\frac{16}{(6n+1)^2}-\frac{24}{(6n+2)^2}-\frac{8}{(6n+3)^2}-\frac{6}{(6n+4)^2}+\frac{1}{(6n+5)^2}
\right).
\end{equation}

\section{Pretty infinite expansions for $\pi$ using products and radicals}

In $1999$, Osler proved the following formula \cite{osler}:
\begin{equation}\label{osler}
\frac{2}{\pi}=\prod_{n=1}^{p} \sqrt{\frac{1}{2}+\frac{1}{2}
\sqrt{\frac{1}{2}+ \frac{1}{2}
\sqrt{\frac{1}{2}+\stackrel{n}{\cdots}+\frac{1}{2}
\sqrt{\frac{1}{2}}}}} \cdot \prod_{n=1}^{\infty}
\frac{2^{p+1}n-1}{2^{p+1}n} \cdot \frac{2^{p+1}n+1}{2^{p+1}n},
\end{equation}
where $n$ inside the square-root means the number of terms in the sum. This formula includes Vi\`{e}te's formula (\ref{viete}) ($p=\infty$); Wallis' product (\ref{wallis}) ($p=0$) and new formulas. For example, with $p=2$ he gets
\begin{equation}\label{ejemnuevaosler}
\frac{2}{\pi}=\sqrt{\frac{1}{2}} \sqrt {\frac{1}{2}+\frac{1}{2}
\sqrt{ \frac{1}{2}}} \cdot \frac{7}{8} \cdot \frac{9}{8} \cdot \frac{15}{16} \cdot \frac{17}{16} \cdot \frac{23}{24} \cdot \frac{25}{24} \cdots.
\end{equation}
In order to reach his formula, Osler used the double angle sine formula $p$ times in a row:
\[
\sin x=2 \cos \frac{x}{2} \sin \frac{x}{2}=2^2 \cos \frac{x}{2} \cos \frac{x}{2^2} \sin \frac{x}{2^2}=
2^p \cos \frac{x}{2} \cos \frac{x}{2^2} \cdots \cos \frac{x}{2^p}
\sin \frac{x}{2^p}
\]
and substituted $t=\frac{x}{2^p}$ in Euler's formula
\[ \frac{\sin t}{t}=\prod_{n=1}^{\infty} \left( 1- \frac{t^2}{\pi^2 n^2} \right)=
\prod_{n=1}^{\infty} \left( \frac{n \pi-t}{n \pi} \cdot \frac{n
\pi+t}{n \pi} \right).
\]
Combining the preceding formulas he got
\[
\frac{\sin x}{x}=\cos \frac{x}{2} \cos \frac{x}{2^2} \cdots \cos \frac{x}{2^p}
\cdot \prod_{n=1}^{\infty} \left( \frac{2^p n \pi-x}{2^p n \pi}
\cdot \frac{2^p n \pi+x}{2^p n \pi} \right).
\]
Finally, he used the elementary trigonometry relations
\[ \cos \frac{x}{2}=\sqrt{\frac{1}{2}+\frac{1}{2} \cos x}, \qquad \cos \frac{x}{2^2}=
\sqrt{\frac{1}{2}+\frac{1}{2}\sqrt{\frac{1}{2}+\frac{1}{2} \cos x}},
\qquad \cdots
\]
and substituted $x=\frac{\pi}{2}$.

\par In $2002$, J. Sondow obtained another pretty expression with products and radicals \cite{sondow}:
\begin{equation}\label{sondow}
\frac{\pi}{2}=\prod_{n=0}^{\infty} \left[ 1^{(-1)^1 {n \choose 0}}
\cdot 2^{(-1)^2 {n \choose 1}} \cdots (n+1)^{(-1)^{n+1} {n \choose
n}} \right]^{\frac{1}{2^{n}}}.
\end{equation}
He started by taking logarithms in Wallis' product (\ref{wallis}):
\[
\ln \frac{\pi}{2}=\sum_{k=0}^{\infty} (-1)^k \ln \frac{k+2}{k+1}.
\]
Taking into account that $\sum_{n \geq k} \frac{{n \choose k}}{2^{n+1}}=1$, we continue like this:
\[
\ln \frac{\pi}{2}=\sum_{k \geq 0} (-1)^k \ln \frac{k+2}{k+1}
\sum_{n \geq k} \frac{{n \choose k}}{2^{n+1}}=\sum_{n=0}^{\infty} {1
\over 2^{n+1}} \sum_{k=0}^n (-1)^k {n \choose k} \ln
\frac{k+2}{k+1}.
\]
Then he substituted $n-1$ for $n$, rewrote the logarithm as $\ln(k+2)-\ln(k+1)$ and the sum on $k$ as the difference of two sums. He replaced $k$ with $k-1$ in the first of these two sums. Then he applied the binomial relation ${n-1 \choose k-1}+{n-1 \choose k}={n \choose k}$ and by exponentiating the expression he got the result.

\section{Very recent results: Series for $1/\pi$, $1/\pi^2$ and $1/\pi^3$}

T. Sato, in a 2002 conference in Japan, showed the series
\begin{equation}\label{sato}
\sum_{n=0}^{\infty} u_n \left( \frac{\sqrt{5}-1}{2} \right)^{12n}
(20n+10-3\sqrt{5}) =\frac{20 \sqrt{3}+9 \sqrt{15}}{6 \pi},
\end{equation}
where $u_n$ are the Ap\'{e}ry numbers, defined by
\[ u_n=\sum_{k=0}^n {n \choose k}^2 {n+k \choose k}^2. \]
Soon afterwards this type of series was found for other types of numbers \cite{chan2}, \cite{yang}. H. H. Chan and Y. Yang independently developed a theory explaining this kind of formula, which also includes as a particular case Ramanujan type series, and they named them Ramanujan-Sato type formulas.

\par A radically different explanation for some Ramanujan type formulas, one that does not use the theory of modular functions, is based on the WZ method, due to Herb Wilf and Doron Zeilberger \cite{petkovek}. It lets one automatically prove identities of the shape
\[ \sum_{n=0}^{\infty} G(n,k)={\rm Constant} \]
in case the function $G(n,k)$ is hypergeometric in both of its variables, that is, if the quotients $G(n+1,k)/G(n,k)$ and $G(n,k+1)/G(n,k)$ are rational functions of $n$ and $k$. The program EKHAD, written by D. Zeilberger, does the hard work of finding a companion function $F(n,k)$ to $G(n,k)$ such that $F(0,k)=0$, and which makes $(F,G)$ a Wilf and Zeilberger (WZ) pair, that is, a pair characterized by the property
\[ F(n+1,k)-F(n,k)=G(n,k+1)-G(n,k). \]
In addition, Zeilberger has proved that if we define $H(n,k)=F(n+1,n+k)+G(n,n+k)$, then we will obtain an identity associated to the first one:
\[ \sum_{n=0}^{\infty} H(n,k)=\sum_{n=0}^{\infty} G(n,k)= {\rm Constant}. \]

\par It is important to emphasize that the WZ method is a symbolic method and not a numerical one like the PSLQ algorithm (see section 9), so the proofs obtained are completely rigorous. The problem is that it is a difficult task to discover suitable functions belonging to a WZ pair.

\par Some identities discovered by the author  \cite{guillera2}, among them
\begin{equation}\label{gui1}
\sum_{n=0}^{\infty} \frac{1}{2^{8n}2^{4k}} \frac{{2n \choose n}^3
{2k \choose k}^2}{{n+k \choose k}^2} (6n+4k+1)=\frac{4}{\pi},
\end{equation}
\begin{equation}\label{gui2}
\sum_{n=0}^{\infty} \frac{(-1)^n}{2^{12n} 2^{8k}} \frac{{2n \choose
n}^5 {2k \choose k}^4}{{n+k \choose k}^4} (20n^2+8n+1+24kn+8k^2+4k)=
\frac{8}{\pi^2},
\end{equation}
are suitable for having an automatic proof with EKHAD. Just to see how the WZ method works, we explain in detail the proof of the identity (\ref{gui1}) based on this method. Given the function
\[ G(n,k)= \frac{1}{2^{8n}2^{4k}} \frac{{2n \choose n}^3 {2k \choose k}^2}{{n+k \choose k}^2} (6n+4k+1), \]
EKHAD gives the companion
\[ F(n,k)=\frac{1}{2^{8n}2^{4k}} \frac{{2n \choose n}^3 {2k \choose k}^2}{{n+k \choose k}^2} 8n \]
and together they form a WZ pair, that is, they satisfy
\[ F(n+1,k)-F(n,k)=G(n,k+1)-G(n,k). \]
We sum over all $n \geq 0$, so that a telescoping cancellation is produced on the first member, and we get
\[  \sum_{n=0}^{\infty} G(n,k) = \sum_{n=0}^{\infty} G(n,k+1). \]
This series converges uniformly, so we can write
\[
\lim_{k \to \infty} \sum_{n=0}^{\infty} G(n,k)=\sum_{n=0}^{\infty}
\lim_{k \to \infty} G(n,k)=\lim_{k \to \infty} G(0,k)=\lim_{k \to
\infty} \frac{{2k \choose k}^2}{2^{4k}} (4k+1) =\frac{4}{\pi}
\]
and we obtain the result
\[  \sum_{n=0}^{\infty} G(n,k)=\frac{4}{\pi}. \]
And, as explained before, we have the double identity
\[ \sum_{n=0}^{\infty} G(n,k)=\sum_{n=0}^{\infty} H(n,k)= \frac{4}{\pi}, \]
where $H(n,k)=F(n+1,n+k)+G(n,n+k).$
Substituting $k=0$, we obtain the Ramanujan series (\ref{rama42n+5})
\begin{equation}\label{ramanujanwz}
\sum_{n=0}^{\infty}  \frac{(2n)!^3}{n!^6} \frac{1}{2^{8n}} (6n+1)=
\frac{1}{4} \sum_{n=0}^{\infty} \frac{(2n)!^3}{n!^6}
\frac{1}{2^{12n}} (42n+5)=\frac{4}{\pi}.
\end{equation}

\par For the function $G(n,k)$ in the series (\ref{gui2}), EKHAD finds its companion
\[
F(n,k)=\frac{(-1)^n}{2^{12n} 2^{8k}} \frac{{2n \choose n}^5
{2k \choose k}^4}{{n+k \choose k}^4} \cdot 8n (2n+4k+1),
\]
and then we can say that the proof of (\ref{gui2}) is hidden inside this pair $(F,G)$.
We have the double identity
\[ \sum_{n=0}^{\infty} G(n,k)=\sum_{n=0}^{\infty} H(n,k)= \frac{8}{\pi^2}. \]
Substituting $k=0$, we obtain the new formulas \cite{guillera1}, \cite{guillera2}, \cite{guilleratesis}
\begin{equation}\label{gui3}
\sum_{n=0}^{\infty} \frac{(2n)!^5}{n!^{10}} \frac{(-1)^n}{2^{12n}}
(20n^2+8n+1)= \frac{1}{16} \sum_{n=0}^{\infty}
\frac{(2n)!^5}{n!^{10}} \frac{(-1)^n}{2^{20n}}
(820n^2+180n+13)=\frac{8}{\pi^2}.
\end{equation}

\par Using the PSLQ algorithm, the author managed to find more series of this type for $1/\pi^2$ \cite{guillera3}, for example,
\begin{equation}\label{gui4}
\sum_{n=0}^{\infty} \frac{(6n)!}{n!^6} \frac{(-1)^n}{2880^{3n}}
(5418n^2+693n+29)={128 \sqrt 5 \over \pi^2},
\end{equation}
whose proof is still unknown. In \cite{guillera4} and \cite{guilleratesis} we found the same series as the result of a conjecture.

\par B. Gourevitch found a series for $1/\pi^3$ \cite{guillera3} using the PSLQ algorithm
\begin{equation}\label{boris}
\sum_{n=0}^{\infty} \frac{(2n)!^7}{n!^{14}} \frac{1}{2^{20n}} (168
n^3+76 n^2+14n+1)={32 \over \pi^3}.
\end{equation}
It has not been proved yet.

\par W. Zudilin proved the following new identities \cite{zudilin} by means of quadratic transformations of the identities (\ref{gui3}):
\begin{equation}\label{zudilin1}
\sum_{n=0}^{\infty} w_n \frac{(4n)!}{n!^2(2n)!} (18n^2-10n-3)
\frac{1}{(2^8 5^2)^n} =\frac{10 \sqrt{5}}{\pi^2},
\end{equation}
\begin{equation}\label{zudilin2}
\sum_{n=0}^{\infty} w_n \frac{(4n)!}{n!^2(2n)!}
(1046529n^2+227104n+16032) \frac{1}{(5^4 41^2)^n}=\frac{5^4 41
\sqrt{41}}{\pi^2},
\end{equation}
where the $w_n$ are numbers defined by
\[ w_n=\sum_{k=0}^{n} {2k \choose k}^3 {2n-2k \choose n-k} 2^{4(n-k)}. \]
With these last formulas we conclude our historical review.

\section*{Acknowledgements} I thank Professors Javier Cilleruelo and Jonathan Sondow for their helpful suggestions and comments. I thank Luis S\'{a}nchez Lajusticia for translating the original Spanish manuscript.

{\it Av. Ces\'{a}reo Alierta, 31 esc. izda {\rm $4^o$}--A, Zaragoza
50008, Spain. \par jguillera@gmail.com}

\end{document}